\documentclass[12pt]{article}

\usepackage{latexsym, amssymb, amsmath, amscd, amsfonts, epsfig, graphicx, colordvi,verbatim,ifpdf}
\usepackage{amsfonts, amsmath, amssymb}
\usepackage{amssymb,amsfonts,amsmath,latexsym,epsfig,cite, psfrag,eepic,color}
\usepackage{amscd,graphics}
\usepackage{latexsym, amssymb,  amsmath,amscd, amsfonts, epsfig, graphicx, colordvi,amsthm}

\usepackage{graphicx}
\usepackage{color}
\usepackage{ifpdf}
\usepackage{fancybox}

\newtheorem{thm}{Theorem}[section]

\newtheorem{coro}[thm]{Corollary}

\def\pf{\noindent{\it Proof.} }
\def\qed{\nopagebreak\hfill{\rule{4pt}{7pt}}
\medbreak}

\numberwithin{equation}{section}

\def\qed{\nopagebreak\hfill{\rule{4pt}{7pt}}
\medbreak}

\newlength{\boxedparwidth}
  {\begin{center} \begin{tabular}{|@{\hspace{.315in}}c@{\hspace{.15in}}|}
                  \hline \\ \begin{minipage}[t]{\boxedparwidth}
                  \setlength{\parindent}{.25in}}%
  {\end{minipage} \\ \\ \hline \end{tabular} \end{center}}

\begin{document}
\begin{center}

{ \large\bf On the Positive Moments of  Ranks of Partitions}
\end{center}

\vskip 5mm

\begin{center}
{  William Y.C. Chen}$^{1}$ ,    {  Kathy Q. Ji}$^{2}$
and {   Erin  Y.Y. Shen}$^{3}$ \vskip 2mm

    $^{1,2,3}$Center for Combinatorics, LPMC-TJKLC\\
   Nankai University, Tianjin 300071, P.R. China\\[6pt]
   $^{1}$Center for Applied Mathematics\\
Tianjin University,  Tianjin 300072, P. R. China\\[6pt]

   \vskip 2mm

    $^1$chen@nankai.edu.cn, $^2$ji@nankai.edu.cn,  $^3$shenyiying@mail.nankai.edu.cn
\end{center}

\vskip 6mm \noindent {\bf Abstract.} By introducing
$k$-marked Durfee symbols, Andrews found a
 combinatorial interpretation of   $2k$-th symmetrized  moment $\eta_{2k}(n)$ of ranks of partitions of $n$ in terms of $(k+1)$-marked Durfee symbols of $n$.      In this paper, we consider the $k$-th symmetrized positive moment $\bar{\eta}_k(n)$ of ranks of partitions of $n$ which is
 defined as the truncated sum over positive ranks of partitions of $n$.   As combintorial interpretations of  $\bar{\eta}_{2k}(n)$ and $\bar{\eta}_{2k-1}(n)$, we show that for fixed $k$ and $i$ with  $1\leq i\leq k+1$,  $\bar{\eta}_{2k-1}(n)$ equals the number of $(k+1)$-marked Durfee symbols of $n$  with the   $i$-th rank  being zero  and $\bar{\eta}_{2k}(n)$ equals the number of $(k+1)$-marked Durfee symbols of $n$ with
  the $i$-th rank being positive. The interpretations of $\bar{\eta}_{2k-1}(n)$ and $\bar{\eta}_{2k}(n)$ also imply
  the interpretation of $\eta_{2k}(n)$ given by Andrews
  since $\eta_{2k}(n)$ equals $\bar{\eta}_{2k-1}(n)$ plus twice of $\bar{\eta}_{2k}(n)$.   Moreover,  we obtain the generating functions of $\bar{\eta}_{2k}(n)$ and $\bar{\eta}_{2k-1}(n)$.



\section{Introduction}

This paper is concerned with a combinatorial study of the symmetrized positive moments of ranks of partitions.
The notion of symmetrized moments was introduced
by Andrews \cite{Andrews-07-a}. The  odd   symmetrized moments  are zero due to the symmetry of ranks. For even symmetrized moments, Andrews found a combinatorial interpretation  by introducing $k$-marked Durfee symbols.
It is natural to
 investigate the combinatorial interpretation of the odd symmetrized moments
 which are   truncated sum over positive ranks of partitions of $n$.  We give  combinatorial interpretations of the even and odd  positive moments   in terms of $k$-marked Durfee symbols, which also lead to the combinatorial interpretation of the even symmetrized moments of ranks given by Andrews.


The rank of a partition $\lambda$ introduced by Dyson \cite{Dyson-1944} is defined as the largest part minus the number of parts.  Let $N(m,n)$ denote the number of partitions of $n$ with rank $m$. The generating function of $N(m,n)$ is given by
\begin{thm}[Dyson-Atkin-Swinnerton-Dyer \cite{Atki54}]\label{gf-r} For fixed integer $m$, we have
\begin{equation}\label{gf-c-e}
\sum_{n=0}^{+\infty}
N(m,n)q^n=\frac{1}{(q;q)_{\infty}}\sum_{n=1}^{+\infty}(-1)^{n-1}
q^{n(3n-1)/2+|m|n}(1-q^n).
\end{equation}
\end{thm}

Recently,  Andrews \cite{Andrews-07-a} introduced the $k$-th symmetrized moment $\eta_{k}(n)$ of ranks  of partitions of $n$ as given by
\begin{align}\label{symRankMom}
\eta_{k}(n)=\sum_{m=-\infty}^{+\infty}{m+\lfloor\frac{k-1}{2}\rfloor \choose k}N(m,n).
\end{align}
 It can be easily seen that for given $k$, $\eta_{k}(n)$ is  a linear combination of the moments $N_j(n)$ of ranks given by Atkin and Garvan \cite{Atkin-Garvan-03}
\begin{align*}
N_j(n)=\sum_{m=-\infty}^{\infty}m^jN(m,n).
\end{align*}
For example,
\[\eta_6(n)=\frac{1}{720}N_6(n)-\frac{1}{144}N_4(n)+\frac{1}{180}N_2(n).\]

In view of the symmetry $N(-m,n)=N(m,n)$,  we have $\eta_{2k+1}(n)=0$.  As for the even symmetrized
moments $\eta_{2k}(n)$, Andrews  gave the following combinatorial interpretation by introducing $k$-marked Durfee symbols. For the definition of $k$-marked Durfee symbols, see Section 2.

\begin{thm}[Andrews \cite{Andrews-07-a}]\label{Andrews} For fixed $k\geq 1$,
  $\eta_{2k}(n)$  is equal to the number of $(k+1)$-marked Durfee symbols of $n$.
\end{thm}

 Andrews  \cite{Andrews-07-a}  proved the above theorem by using the $k$-fold generalization of Watson's $q$-analog of Whipple's theorem. Ji \cite{Ji-2011} gave a combinatorial proof of Theorem \ref{Andrews} by  establishing a map from  $k$-marked Durfee symbols to ordinary partitions.  Kursungoz \cite{Kursungoz-2011} provided another proof of Theorem \ref{Andrews}  by using an alternative representation of $k$-marked Durfee symbols.

In this paper,  we introduce the $k$-th symmetrized positive moment $\bar{\eta}_{k}(n)$ of ranks as given by
\begin{align*}
\overline{\eta}_k(n)=\sum_{m=1}^{\infty}
{m+\lfloor\frac{k-1}{2}\rfloor\choose k}N(m,n),
\end{align*}
or equivalently,
\begin{equation}
\overline{\eta}_{2k-1}(n)=\sum_{m=1}^{\infty}
{m+k-1\choose 2k-1}N(m,n)
\end{equation}
and
\begin{equation}
\overline{\eta}_{2k}(n)=\sum_{m=1}^{\infty}
{m+k-1\choose 2k}N(m,n).
\end{equation}
 Furthermore, it is easy to see that for given $k$, $\bar{\eta}_k(n)$ is a linear combination of the positive moments $\overline{N}_j(n)$ of ranks introduced by  Andrews, Chan and Kim \cite{Andrews-Chan-Kim} as given by
\begin{align*}
\overline{N}_j(n)=\sum_{m=1}^{\infty}m^jN(m,n).
\end{align*}
For example,
\begin{align*}
\bar{\eta}_4(n)&=
\frac{1}{24}\overline N_4(n)-\frac{1}{12}\overline N_3(n)-\frac{1}{24}\overline N_2(n)
+\frac{1}{12}\overline N_1(n),\\[5pt]
\bar{\eta}_5(n)&=\frac{1}{120}\overline{N}_5(n)-\frac{1}{24}\overline{N}_3(n)
+\frac{1}{30}\overline{N}_1(n).
\end{align*}
By the symmetry $N(-m,n)=N(m,n)$, it is readily seen that
 \begin{equation}\label{re}
 \eta _{2k}(n)=2 \overline{\eta}_{2k}(n)+ \overline{\eta} _{2k-1}(n).
\end{equation}

 The main objective of this paper is to give combintorial interpretations of  $\bar{\eta}_{2k}(n)$ and $\bar{\eta}_{2k-1}(n)$. We show that for given $k$ and $i$ with  $1\leq i\leq k+1$,  $\bar{\eta}_{2k-1}(n)$ equals the number of $(k+1)$-marked Durfee symbols of $n$  with the   $i$-th rank  being zero  and $\bar{\eta}_{2k}(n)$ equals the number of $(k+1)$-marked Durfee symbols of $n$ with
  the $i$-th rank being positive. It should be noted that
  $\bar{\eta}_{2k-1}(n)$ and $\bar{\eta}_{2k}(n)$ are independent of $i$
  since the ranks of $k$-marked Durfee symbols are symmetric, see Andrews \cite[Corollary 12]{Andrews-07-a}.

   With the aid of  Theorem \ref{main-1} and Theorem \ref{main-2} together with the   generating function \eqref{gf-c-e} of $N(m,n)$,  we  obtain the generating functions of  $\bar{\eta}_{2k}(n)$ and  $\bar{\eta}_{2k-1}(n)$.

\section{Combinatorial interpretations}

 In this section, we   give  combinatorial interpretations of $\bar{\eta}_{2k-1}(n)$ and $\bar{\eta}_{2k}(n)$ in terms of the $k$-marked Durfee symbols.
For a partition $\lambda$, we write $\lambda=(\lambda_1, \lambda_2, \ldots, \lambda_s)$, so that
$\lambda_1$ is the largest part and $\lambda_s$ is the smallest part of $\lambda$. Recall that a $k$-marked Durfee symbol  of $n$ introduced by Andrews \cite{Andrews-07-a} is a two-line array composed of $k$ pairs $(\alpha^i,\beta^i)$ of
partitions along with a positive integer $D$ which
 is represented in the following form:
  \[\tau=\left(\begin{array}{cccc}
\alpha^k,&\alpha^{k-1},&\ldots,&\alpha^1\\[2pt]
\beta^k,&\beta^{k-1},&\ldots,&\beta^1
\end{array}
\right)_D,\] where  the partitions $\alpha^i$ and $\beta^i$  satisfy
the following four conditions:
\begin{itemize}

\item[{\rm (1)}]The partitions $\alpha^i$ ($1\leq i<k$) are nonempty, while $\alpha^k$ and  $\beta^i$ ($1\leq i\leq k$) are allowed to be empty{\rm ;}

\item[{\rm (2)}] $\beta^{i-1}_1\leq \alpha^{i-1}_1 \leq
\min\{\alpha^i_s,\beta^{i}_{s}\}$ for $2\leq i\leq k${\rm ;}

\item[{\rm (3)}]$\alpha_1^k$, $\beta_1^k \leq D${\rm ;}

\item[{\rm (4)}]$\sum_{i=1}^{k}(|\alpha^i|+|\beta^i|)+D^2=n$.
\end{itemize}

Let
\[\tau=\left(\begin{array}{cccc}
\alpha^k,&\alpha^{k-1},&\ldots,&\alpha^1\\[2pt]
\beta^k,&\beta^{k-1},&\ldots,&\beta^1
\end{array}
\right)_D\]
be a $k$-marked Durfee symbol. The pair
$(\alpha^i,\beta^i)$ of partitions
is called the $i$-th vector of $\tau$. Andrews defined
 the $i$-th rank $\rho_i(\tau)$ of $\tau$ as follows
\[\rho_i(\tau)=\left\{
\begin{array}{ll}
\ell(\alpha^i)-\ell(\beta^i)-1, \ \ &\text{ for }\ 1\leq i<k,\\[5pt]
\ell(\alpha^k)-\ell(\beta^k). \ \ &\text{ for }\ i=k.
\end{array}\right.
\]

For example, consider  the following   $3$-marked Durfee symbol $\tau$.
\[\tau=\left(\begin{array}{cccccccc}
\overbrace{5_3,4_3}^{\alpha^3},&
\overbrace{4_2,3_2,3_2,2_2}^{\alpha^2},&
\overbrace{2_1}^{\alpha^1}\\[3pt]
\underbrace{4_3}_{\beta^3},&
\underbrace{3_2,2_2,2_2}_{\beta^2},&
\underbrace{2_1,2_1}_{\beta^1}
\end{array}\right)_5.\]
We have  $\rho_1(\tau)=-2,\,\rho_2(\tau)=0,$ and $\rho_3(\tau)=1$.

For odd symmetrized moments $\bar{\eta}_{2k-1}(n)$, we have the following combinatorial interpretation.

\begin{thm}\label{main-1} For fixed positive integers  $k$ and $i$ with $1\leq i\leq k+1$,
 $\bar{\eta}_{2k-1}(n)$ is equal to the number of $(k+1)$-marked Durfee symbols of $n$ with the $i$-th rank equal to zero.
\end{thm}

For the even case, we
have the following interpretation.

\begin{thm}\label{main-2}For fixed positive integers  $k$ and $i$ with $1\leq i\leq k+1$,
 $\bar{\eta}_{2k}(n)$  is equal to the number of $(k+1)$-marked Durfee symbols of $n$ with the $i$-th rank being positive.
\end{thm}

The proofs of the above two interpretations are based on the following partition identity given by Ji  \cite{Ji-2011}.
We shall adopt the notation $ {D}_k(m_1,m_2,\ldots,m_k;n)$ as used by
Andrews \cite{Andrews-07-a} to denote the   number of
$k$-marked Durfee symbols of $n$ with $i$-th rank equal to
$m_i$.

 \begin{thm}\label{r-thm} Given $k\geq 2$ and $n\geq 1$, we have
\begin{align}\label{r-thm-e}
 {D}_{k}(m_1,m_2,\ldots,m_k;n)
=\sum_{t_1,\ldots,t_{k-1}=0}^{\infty}
N\left(\sum_{i=1}^{k}|m_i|+2\sum_{i=1}^{k-1}t_i+k-1,n\right).
\end{align}
\end{thm}

To prove the above two interpretations, we also need
 the following symmetric property given by Andrews \cite{Andrews-07-a}.  Boulet and Kursungoz \cite{Boulet-2011} found a combinatorial proof of this fact.

\begin{thm} \label{r-sym-2}For $k\geq 2$ and $n\geq 1$,
$ {D}_k(m_1,\ldots,m_k;n)$ is symmetric in $m_1,\ m_2,\ldots,m_k$.
\end{thm}

We are now in a position to prove Theorem \ref{main-1} and Theorem \ref{main-2} with the aid of  Theorem \ref{r-thm} and Theorem \ref{r-sym-2}.

\noindent {\it Proof of Theorem \ref{main-1}.}
By Theorem \ref{r-sym-2}, it suffices to show that
\begin{equation}\label{step-1}
\sum_{m_2,m_3,\ldots,m_{k+1}=-\infty}^{\infty}
 {D}_{k+1}(0,m_2,m_3,\ldots,m_{k+1};n)=\bar{\eta}_{2k-1}(n).
\end{equation}
Using Theorem \ref{r-thm}, we get
\begin{align}\label{main-1-e}
&\sum_{m_2,m_3,\ldots,m_{k+1}=-\infty}^{\infty}
 {D}_{k+1}(0,m_2,m_3,\ldots,m_{k+1};n)\nonumber \\
&\hskip 2cm =\sum_{m_2,m_3,\ldots,m_{k+1}=-\infty}^{\infty}
\sum_{t_1,\ldots,t_{k}=0}^{\infty}
N\left(\sum_{i=2}^{k+1}|m_i|+2\sum_{i=1}^{k}t_i+k,n\right).
\end{align}

Given $k$ and $n$, let $c_k(n)$ denote the number  of
 integer solutions to the equation
\[|m_2|+\cdots+|m_{k+1}|+2t_1+\cdots+2t_k=n,\]
where the variables   $m_i$ are integers and the variables $t_i$ are nonnegative integers. It is easy to see that the generating function of $c_k(n)$ is equal to
\begin{align}
\sum_{n=0}^{\infty}c_k(n)q^n&=(1+2q+2q^2+2q^3+\cdots)^{k}
(1+q^2+q^4+q^6+\cdots)^{k}\nonumber\\[3pt]
&=\left(\frac{1+q}{1-q}\right)^k\left(\frac{1}{1-q^2}\right)^k
\nonumber\\[3pt]
&=\frac{1}{(1-q)^{2k}}\nonumber\\[3pt]
&=\sum_{n=0}^{\infty}{n+2k-1 \choose 2k-1}q^n. \label{gf-ck}
\end{align}
 Equating the coefficients of $q^n$ on the both sides of \eqref{gf-ck}, we get
\begin{align*}
c_k(n)={n+2k-1 \choose 2k-1},
\end{align*}
that is,
\[c_k(m-k)={m+k-1 \choose 2k-1}.\]
Thus, \eqref{main-1-e} can be written as
\begin{align*}
&\sum_{m_2,m_3,\ldots,m_{k+1}=-\infty}^{\infty}
 {D}_{k+1}(0,m_2,m_3,\ldots,m_{k+1};n)\nonumber \\
&\hskip 2cm = \sum_{m=1}^{\infty}
{m+k-1 \choose 2k-1}N(m,n)
\end{align*}
which is equal to $\bar{\eta}_{2k-1}(n).$  This completes the proof.  \qed

\medskip

\noindent {\it Proof of Theorem \ref{main-2}.}  Similarly, by Theorem \ref{r-sym-2},
it is enough  to show that
\begin{equation}\label{main-2-1}
\sum_{\stackrel{m_1>0}{m_2,m_3,\ldots,m_{k+1}=-\infty}}^{\infty}
 {D}_{k+1}(m_1,m_2,\ldots,m_{k+1};n)=\bar{\eta}_{2k}(n).
\end{equation}
Using Theorem \ref{r-thm}, we get
\begin{align}\label{last-main-2}
&\sum_{\stackrel{m_1>0}{m_2,m_3,\ldots,m_{k+1}=-\infty}}^{\infty}
 {D}_{k+1}(m_1,m_2,\ldots,m_{k+1};n)\nonumber \\[3pt]
&\hskip 1cm \qquad =\sum_{\stackrel{m_1>0}{m_2,m_3,\ldots,m_{k+1}=-\infty}}^{\infty}\sum_{t_1,\ldots,t_{k}=0}^{\infty}
N\left(m_1+\sum_{i=2}^{k+1}|m_i|+2\sum_{i=1}^{k}t_i+k,n\right).
\end{align}
Given $k$ and $n$, let $\bar{c}_k(n)$ denote the number of integer solutions to the equation
\[m_1+|m_2|+\cdots+|m_{k+1}|+2t_1+\cdots+2t_k=n,\]
where the variable $m_1$ is a positive integer, the variables $m_i$ ($2\leq i\leq k+1$) are   integers and the variables $t_i$ are  nonnegative integers. An easy computation shows that
\begin{align}
 \sum_{n=0}^{\infty}\bar{c}_k(n)q^n =
 \frac{q}{(1-q)^{2k+1}},
\end{align}
so that
\begin{align*}
\bar{c}_k(n)={n+2k-1 \choose 2k }.
\end{align*}
We write
\[\bar{c}_k(m-k)={m+k-1 \choose 2k}.\]
It follows that
\begin{align*}
&\sum_{\stackrel{m_1>0}{m_2,m_3,\ldots,m_{k+1}=-\infty}}^{\infty}
 {D}_{k+1}(m_1,m_2,\ldots,m_{k+1};n)\nonumber \\
&\hskip 2cm = \sum_{m=1}^{\infty}
{m+k-1 \choose 2k}N(m,n),
\end{align*}
which equals $\bar{\eta}_{2k}(n),$  as required.   \qed

Note that   the number ${D}_k(m_1,\ldots,
m_k;n)$ has the mirror symmetry with
respect to each $m_i$, that is, for $1\leq i\leq k$, we have
\[
 {D}_k(m_1,\ldots,m_i,\ldots,
m_k;n)= {D}_k(m_1,\ldots,-m_i,\ldots, m_k;n).
\]
  Using this mirror symmetry, Theorem \ref{main-2} can be restated  as follows.

\begin{thm}\label{main-3}For fixed positive integers  $k$ and $i$ with $1\leq i\leq k+1$,
 $\bar{\eta}_{2k}(n)$ is also equal to the number of $(k+1)$-marked Durfee symbols of $n$ with the $i$-th rank being negative.
\end{thm}

 \begin{table}[h]
 \[
\begin{array}{c|c|c|c }
&\overline{\eta}_1(5)&\overline{\eta}_2(5)&\overline{\eta}_2(5)\\[5pt]
\hline
 &\left(\begin{array}{llll}
 1_2 &1_2 &1_2 &1_1\\
&&&\end{array}\right)_1&
\left(\begin{array}{llll}
 1_1 &1_1 &1_1 &1_1\\
 &&&\end{array}\right)_1&
\left(\begin{array}{llll}
  1_1&& \\
1_1 &1_1 &1_1\end{array}\right)_1\\[20pt]
& \left(\begin{array}{llll}
 1_2 &1_1 &1_1 \\
1_1&\end{array}\right)_1&\left(\begin{array}{llll}
 1_2 &1_1 &1_1&1_1 \\
&\end{array}\right)_1&\left(\begin{array}{llll}
 1_2 &1_1 \\
1_1&1_1 \end{array}\right)_1\\[20pt]
&\left(\begin{array}{llll}
1_2 &1_2&1_1\\
 1_2 &\end{array}\right)_1&\left(\begin{array}{llll}
1_2 &1_2&1_1&1_1\\
 &\end{array}\right)_1&\left(\begin{array}{llll}
1_2 &1_2&1_1 \\
1_1\end{array}\right)_1\\[20pt]
& \left(\begin{array}{llll}
1_1&\\
 1_2 &1_2&1_2\end{array}\right)_1&\left(\begin{array}{llll}
1_1&1_1&1_1\\
1_1\end{array}\right)_1&\left(\begin{array}{llll}
1_1&1_1\\
1_1&1_1\end{array}\right)_1\\[20pt]
&\left(\begin{array}{llll}
1_1 &1_1\\
1_2& 1_1\end{array}\right)_1&\left(\begin{array}{llll}
1_1 &1_1&1_1\\
1_2 \end{array}\right)_1&\left(\begin{array}{llll}
1_1 \\
1_2&1_1&1_1 \end{array}\right)_1\\[20pt]
&\left(\begin{array}{llll}
 1_2 &1_1 \\
1_2 &1_2\end{array}\right)_1&\left(\begin{array}{llll}
 1_2 &1_1&1_1 \\
1_2 \end{array}\right)_1&\left(\begin{array}{llll}
 1_2 &1_1 \\
1_2 &1_1\end{array}\right)_1\\[20pt]
&\left(\begin{array}{llll}
1_1\\
\ \end{array}\right)_2 &\left(\begin{array}{llll}
1_1&1_1\\
 1_2 &1_2\end{array}\right)_1&\left(\begin{array}{llll}
1_1 \\
 1_2 &1_2&1_1\end{array}\right)_1
 \end{array}
 \]
 \caption{$2$-Marked Durfee Symbols of $5$.}
\end{table}
For example,  for   $n=5$, $k=1$ and $i=1$,  there are twenty-one $2$-marked Durfee symbols of $5$ as listed in Table 2.1. The first column  in Table 2.1 gives seven  $2$-marked Durfee symbols $\tau$ with  $\rho_1(\tau)=0$,  the second column  contains seven $2$-marked Durfee symbols  $\tau$ with $\rho_1(\tau)>0$ and the third column contains seven $2$-marked Durfee symbols $\tau$ with  $\rho_1(\tau)<0$.  It can be verified that $\overline{\eta}_{1}(5)=7$, $\overline{\eta}_{2}(5)=7$ and $\eta_{2}(5)=\overline{\eta}_{1}(5)+2\overline{\eta}_{2}(5)=21.$

\section{The generating functions of $\bar{\eta}_{2k-1}(n)$ and $\bar{\eta}_{2k}(n)$}

In this section, we obtain the generating functions of $\bar{\eta}_{2k-1}(n)$ and $\bar{\eta}_{2k}(n)$ with the aid of Theorem \ref{main-1} and Theorem \ref{main-2}. In doing so, we use the generating function  of $N(m,n)$ to derive the  generating functions  of  $D_{k+1}(0,m_2,\ldots,m_{k+1};n)$ and $D_{k+1}(m_1,m_2,\ldots,m_{k+1};n)$ ($m_1>0$).

\begin{thm}\label{GF-thm-1}
For $k \geq 1$, we have
\begin{align}\label{gf-thm-e}
&\sum_{m_2,\ldots,m_{k+1}=-\infty}^{\infty}
\sum_{n=0}^{\infty}D_{k+1}(0,m_2,\ldots,m_{k+1};n)
x_1^{m_2}\cdots x_{k}^{m_{k+1}}q^n \nonumber \\[3pt]
&\quad \quad =\frac{1}{(q;q)_\infty}\sum_{n=1}^{\infty}(-1)^{n-1}
q^{n(3n-1)/2+kn}
\frac{(1-q^n)}{\prod_{j=1}^{k}(1-x_jq^n)(1-x_j^{-1}q^n)}.
\end{align}
\end{thm}

\pf  Let
\[G_{k}(x_1,\ldots,x_{k};q)
=\sum_{m_2,\ldots,m_{k+1}=-\infty}^{\infty}
\sum_{n=0}^{\infty}D_{k+1}(0,m_2,\ldots,m_{k+1};n)x_1^{m_2}\cdots x_{k}^{m_{k+1}}q^n.\]

By Theorem \ref{r-thm}, we have
\begin{align}
&G_{k}(x_1,\ldots,x_{k};q) \nonumber \\[3pt]
&\quad \quad=\sum_{m_2,\ldots,m_{k+1}=-\infty}^{\infty}
\sum_{t_1,\ldots,t_{k}=0}^{\infty}x_1^{m_2}\cdots x_{k}^{m_{k+1}}\sum_{n=0}^{\infty}
N\left(\sum_{i=2}^{k+1}|m_i|+2\sum_{i=1}^{k}t_i+k,n\right)
q^n. \label{gf-temp}
\end{align}
Using \eqref{gf-c-e} with $m$ replaced by $\sum_{i=2}^{k+1}|m_i|+2\sum_{i=1}^{k}t_i+k$, we get
\begin{align*}\sum_{n=0}^{\infty}
&N\left(\sum_{i=2}^{k+1}|m_i|+2\sum_{i=1}^{k}t_i+k,n\right)q^n\\[3pt]
&\quad =\frac{1}{(q;q)_\infty}\sum_{n=1}^{\infty}(-1)^{n-1}
q^{n(3n-1)/2+n(\sum_{i=2}^{k+1}|m_i|+2\sum_{i=1}^{k}t_i+k)}
(1-q^n).
\end{align*}
Therefore \eqref{gf-temp} becomes
\begin{align}\label{temp}
G_{k}(x_1,\ldots,x_{k};q)&=\sum_{m_2,\ldots,m_{k+1}=-\infty}^{\infty}
\sum_{t_1,\ldots,t_{k}=0}^{\infty}x_1^{m_2}\cdots x_{k}^{m_{k+1}} \nonumber\\[3pt]
& \hskip 2cm \times \frac{1}{(q;q)_\infty}\sum_{n=1}^{\infty}(-1)^{n-1}
q^{n(3n-1)/2+n(\sum_{i=2}^{k+1}|m_i|+2\sum_{i=1}^{k}t_i+k)}
(1-q^n).
\end{align}
Write (\ref{temp}) in the following form
\begin{align}
G_{k}(x_1,\ldots,x_{k};q)
& =\frac{1}{(q;q)_\infty}\sum_{n=1}^{\infty}(-1)^{n-1}
q^{n(3n-1)/2+kn}(1-q^n) \nonumber\\[3pt]
&\hskip 1cm \times \sum_{m_2,\ldots,m_{k+1}=-\infty}^{\infty}
\sum_{t_1,\ldots,t_{k}=0}^{\infty}x_1^{m_2}\cdots x_{k}^{m_{k+1}}
q^{n(\sum_{i=2}^{k+1}|m_i|+2\sum_{i=1}^{k}t_i)}.\label{gf-t-2}
\end{align}
Notice that
\begin{align}\label{temp-1}
&\sum_{a=-\infty}^{{+\infty}}\sum_{b=0}^{{+\infty}}x^{a}q^{n(|a|+2b)}
=\frac{1}{(1-xq^n)(1-x^{-1}q^n)}.
\end{align}
Applying the above formula repeatedly to \eqref{gf-t-2}, we deduce that
\[G_{k}(x_1,\ldots,x_{k};q)=\frac{1}{(q;q)_\infty}\sum_{n=1}^{\infty}(-1)^{n-1}
q^{n(3n-1)/2+kn}
\frac{(1-q^n)}{\prod_{j=1}^{k}(1-x_jq^n)(1-x_j^{-1}q^n)},
\]
as required.   \qed

Setting $x_j=1$ for $1 \leq j\leq k$ in Theorem
\ref{GF-thm-1} and using Theorem \ref{main-1},  we obtain the following generating function of $\bar{\eta}_{2k-1}(n)$.

\begin{coro}
For $k\geq 1$, we have
\begin{align}\label{GF-eq-1}
\sum_{n=1}^{\infty}\bar{\eta}_{2k-1}(n)q^n
=\frac{1}{(q;q)_\infty}\sum_{n=1}^{\infty}(-1)^{n-1}
q^{n(3n-1)/2+kn}\frac{1}{(1-q^n)^{2k-1}}.
\end{align}
\end{coro}

Taking $k=1$ in \eqref{GF-eq-1} and observing that $\bar{\eta}_1(n)=\overline{N}_1(n)$,  we are led to the
generating function   for $\overline{N}_1(n)$ as given  by
Andrews, Chan and Kim in \cite[Theorem 1]{Andrews-Chan-Kim}.

The following generating function can be shown by using the same reasoning as in the proof of Theorem \ref{GF-thm-1}.

\begin{thm}\label{GF-thm-2}
For $k \geq 1$, we have
\begin{align}
&\sum_{\stackrel{m_1>0}{m_2,\ldots,m_{k+1}=-\infty}}^{\infty}
\sum_{n=1}^{\infty}D_{k+1}(m_1,m_2,\ldots,m_{k+1};n)
x_1^{m_1}\cdots x_{k+1}^{m_{k+1}}q^n \nonumber \\[3pt]
&\quad =\frac{1}{(q;q)_\infty}\sum_{n=1}^{\infty}(-1)^{n-1}
q^{n(3n+1)/2+kn}\frac{x_1(1-q^n)}
{(1-x_1q^n)\prod_{j=2}^{k+1}(1-x_jq^n)(1-x_j^{-1}q^n)}.
\end{align}
\end{thm}

Setting   $x_j=1$ for $1 \leq j\leq k+1$ in Theorem
\ref{GF-thm-2} and using Theorem \ref{main-2},
we arrive at the following generating function of $\bar{\eta}_{2k}(n)$.

\begin{coro}
For $k\geq 1$, we have
\begin{align}\label{GF-eq-2}
\sum_{n=1}^{\infty}\bar{\eta}_{2k}(n)q^n
=\frac{1}{(q;q)_\infty}\sum_{n=1}^{\infty}(-1)^{n-1}
q^{n(3n+1)/2+kn}\frac{1}{(1-q^n)^{2k}}.
\end{align}
\end{coro}

\noindent{\bf Acknowledgments.} This work was supported by the 973 Project and the National Science Foundation of China.

\end{document}